\documentstyle[11pt]{article}
\def\Bbb R{{\rm \bf R}}
\def\proclaim#1{\vskip2mm{\bf #1}\em}
\def\endproclaim{\em \vskip2mm}
\def\tag#1{\eqno(#1)}
\def\gathered{\begin{array}{c}}
\def\endgathered{\end{array}}
\def\text{\mbox}

\begin{document}

\title {Extraction formulae for an inverse boundary value problem
for the equation $\nabla\cdot(\sigma-i\omega\epsilon)\nabla u=0$}
\author{Masaru IKEHATA\footnote{
Department of Mathematics,
Faculty of Engineering,
Gunma University, Kiryu 376-8515, JAPAN}}
\maketitle

\begin{abstract}
We consider an inverse boundary value problem for the equation
$\nabla\cdot(\sigma-i\omega\epsilon)\nabla u=0$ in a given bounded domain $\Omega$ at a fixed $\omega>0$.
$\sigma$ and $\epsilon$ denote the conductivity and permittivity of the material forming $\Omega$, respectively.
We give some formulae for extracting information about the location of the discontinuity
surface of $(\sigma,\epsilon)$ from the Dirichlet-to-Neumann map.
In order to obtain results we make use of two methods.
The first is the enclosure method which is based on a new role of the exponentially growing solutions of the equation for 
the background material.  The second is a generalization of the enclosure method
based on a new role of Mittag-Leffler's function.


\end{abstract}


\section{Introduction}

Inject an alternating electric current $j$ across the boundary of a given body $\Omega$.
The resulting voltage potential $u$ inside the body satisfies the equation
$$\displaystyle
\nabla\cdot(\sigma-i\omega\epsilon)\nabla u=0\,\,\text{in}\,\Omega
\tag {1.1}
$$
and the boundary condition 
$$\displaystyle
j=(\sigma-i\omega\epsilon)\nabla u\cdot\nu\vert_{\partial\Omega}.
$$
Here $\sigma=\sigma(x)$ and $\epsilon=\epsilon(x)$ denote the conductivity and permittivity of the body, respectively;
$\omega>0$ denotes the frequency; $\nu$ denotes the unit outward normal vector field to $\partial\Omega$.

This equation can be deduced as an approximation of the system of time-harmonic Maxwell's equations (see Appendix 2 of \cite{SCI})
under the assumption that the magnetic permeability of the body is very small.

In this paper we assume that $\sigma$ and $\epsilon$ on an open set $D$ of $\Omega$ differ from the known constant,
isotropic background conductivity $\sigma_0(x)\equiv\sigma_0$ and permittivity $\epsilon_0(x)\equiv \epsilon_0$;
$D$ is considered an unknown inclusion embedded in $\Omega$.

Briefly, we are interested in the problem of drawing a picture of  $D$ by means of the
observation data. In this paper the observation data means infinitely many pairs  $(u\vert_{\partial D}, j)$ of
the solutions of equation (1.1). This is an idealized formulation of the electrical impedance
tomography. In my opinion, the problem is divided into two parts.

(1) The first problem is to find a formula that extracts useful information about the location
of  $D$ from the observation data without error.

(2) The second problem is that of regularizing the formula: that is, howto modify the formula
when the observation data contain error.

Of course, for practical application, we have to consider the more serious problem: how
to obtain the data needed for the regularized formula from the experimental data.

All the problems mentioned above are important. However, it should be emphasized that
without finding the solution to the first problem nothing can be achieved. In this paper, we
consider the first problem and give two formulae that yield an estimation of $D$ from above.

\subsection{Description of the problem}

Let us formulate our problem more precisely. We consider $\Omega$ a bounded connected open subset
of $\Bbb R^n$, $n = 2,3$ with Lipschitz boundary. In what follows, unless otherwise stated, we assume
that $\sigma$, $\epsilon$ satisfy (A):

$$
\left\{
\begin{array}{c}
\displaystyle
\text{$\sigma$ and $\epsilon$ are $n\times n$ real symmetric matrix-valued functions on $\Omega$};\\
\displaystyle
\text{all components of  $\sigma$ and $\epsilon$ are essentially bounded functions on $\Omega$};\\
\displaystyle
\text{$\sigma$ is non-negative and $\epsilon$ is uniformly positive definite in $\Omega$.}
\end{array}
\right.
\tag {A}
$$

Using the Lax-Milgram theorem and a fact in the spectral theory for the Hermitian operator
in the Hilbert space \cite{M}, we know that, given $f\in H^{1/2}(\partial D)$ there exists the unique weak
solution $u\in H^1(\Omega)$ of the Dirichlet problem
$$\begin{array}{c}
\displaystyle
\nabla\cdot(\sigma-i\omega\epsilon)\nabla u=0\,\,\text{in}\,\Omega,\\
\\
\displaystyle
u=f\,\,\text{on}\,\partial\Omega.
\end{array}
$$

Note that in order to ensure the unique solvability of this boundary value problem it suffices
to assume that one of  $\sigma$ or $\epsilon$ is uniformly positive definite in $\Omega$. 
Define the bounded linear
functional $\Lambda_{\sigma,\epsilon}f$ on $H^{1/2}(\partial\Omega)$ by the formula
$$\displaystyle
<\Lambda_{\sigma,\epsilon}f,g>
=\int_{\Omega}(\sigma-i\omega\epsilon)\nabla u\cdot\nabla vdx
$$
where $g$ is an arbitrary element in $H^{1/2}(\partial\Omega)$ and $v\in H^1(\Omega)$ with $v=g$ on $\partial\Omega$.

Let $D$ be an open subset of $\Omega$ such that $\overline {D}\subset\Omega$.
Assume that $\sigma$, $\epsilon$ take the form

$$\displaystyle
\sigma(x)=\left\{
\begin{array}{rl}
\displaystyle
\sigma_0, & \,\mbox{if $x\in\,\Omega\setminus D$,}\\
\\
\displaystyle
\sigma_0+\alpha (x), & \,\mbox{if $x\in D$;}
\end{array}
\right.
\tag {1.2}
$$
$$\displaystyle
\epsilon(x)=\left\{
\begin{array}{rl}
\displaystyle
\epsilon_0, & \,\mbox{if $x\in\,\Omega\setminus D$,}\\
\\
\displaystyle
\epsilon_0+\beta (x), & \,\mbox{if $x\in D$}
\end{array}
\right.
\tag {1.3}
$$
where both $\sigma_0$ and $\epsilon_0$ are known constants satisfying
$$
\displaystyle
\sigma_0\ge 0;
\tag {1.4}
$$
$$
\displaystyle
\epsilon_0>0.
\tag {1.5}
$$
We assume that both $\alpha(x)$ and $\beta(x)$ together with $D$ are unknown and that
$(\sigma,\epsilon)$ has some kind of discontinuity across  $\partial D$.

{\bf\noindent Problem.}
Find a formula that extracts an information about the location of $D$ from $\Lambda_{\sigma,\epsilon}$.
We call such a formula an extraction formula of the information. In \cite{IE} we considered the
case when $\omega=0$ and gave an extraction formula of the convex hull of $D$; the method predicts
when a plane ($n = 3$), line ($n = 2$) with a given normal vector descending from $\partial\Omega$
hits $\partial D$.

In this paper we give a remark about the applicability of the method. Note that we do not
assume any regularity for either $\alpha$ or $\beta$.

\subsection{A reduction to the case $\sigma_0=1$ and $\epsilon_0=0$}

In this subsection we describe a simple reduction argument.
For $\sigma$ and $\epsilon$ given by (1.2) and (1.3), respectively define

$$\displaystyle
\tilde{\sigma}
=\frac{\sigma_0\sigma+\omega^2\epsilon_0\epsilon}
{\sigma_0^2+\omega^2\epsilon_0^2};
\tag {1.6}
$$
$$\displaystyle
\tilde{\epsilon}
=\frac{\sigma_0\epsilon-\epsilon_0\sigma}
{\sigma_0^2+\omega^2\epsilon_0^2}.
\tag {1.7}
$$
Then we have
$$\displaystyle
\sigma-i\omega\epsilon
=(\sigma_0-i\omega\epsilon_0)
(\tilde{\sigma}-i\omega\tilde{\epsilon}).
\tag {1.8}
$$
Note that $\tilde{\sigma}(x)=1$ and $\tilde{\epsilon}(x)=0$ for
$x\in\Omega\setminus D$.  From (A), (1.4)-(1.6) one knows that $\tilde{\sigma}$ is
uniformly positive definite in $\Omega$.  Then $\Lambda_{\tilde{\sigma},\tilde{\epsilon}}$
is still well defined and from (1.8) one has
$$\displaystyle
\Lambda_{\sigma,\epsilon}
=(\sigma_0-i\omega\epsilon_0)\Lambda_{\tilde{\sigma},\tilde{\epsilon}}.
\tag {1.9}
$$
Therefore, knowing $\Lambda_{\sigma,\epsilon}$ is equivalent to knowing $\Lambda_{\tilde{\sigma},\tilde{\epsilon}}$
through the relationship (1.9).
Moreover, from (1.6) and (1.7) we have
$$\displaystyle
\left(
\begin{array}{c}
\displaystyle
\tilde{\sigma}-1\\
\\
\displaystyle
\tilde{\epsilon}
\end{array}
\right)
=\frac{1}{\sigma_0^2+\omega^2\epsilon_0^2}
\left(
\begin{array}{lr}
\displaystyle
\sigma_0 & \omega^2\epsilon_0\\
\\
\displaystyle
-\epsilon_0 & \sigma_0
\end{array}
\right)
\left(
\begin{array}{c}
\displaystyle
\sigma-\sigma_0\\
\\
\displaystyle
\epsilon-\epsilon_0
\end{array}
\right).
$$
This implies that $(\sigma, \epsilon)$ has a jump from $(\sigma_0, \epsilon_0)$
if and only if $(\tilde{\sigma}, \tilde{\epsilon})$ has a jump from $(1, 0)$.
In particular, if $\sigma_0=0$, then one has
$$\displaystyle
\tilde{\sigma}-1=\frac{\epsilon-\epsilon_0}{\epsilon_0}.
$$
This means the jump of $\tilde{\sigma}$
from $1$ is proportional to the jump of $\epsilon$ from $\epsilon_0$. However,
we do not want to exclude the case when $\sigma_0>0$.
Hereafter we consider the reduced case unless
otherwise stated and therefore one may assume that $\sigma$ is uniformly positive definite in $\Omega$;

$$\displaystyle
\sigma(x)=\left\{
\begin{array}{rl}
\displaystyle
1, & \,\mbox{if $x\in\,\Omega\setminus D$,}\\
\\
\displaystyle
1+a (x), & \,\mbox{if $x\in D$;}
\end{array}
\right.
\tag {1.10}
$$
$$\displaystyle
\epsilon(x)=\left\{
\begin{array}{rl}
\displaystyle
0, & \,\mbox{if $x\in\,\Omega\setminus D$,}\\
\\
\displaystyle
b (x), & \,\mbox{if $x\in D$.}
\end{array}
\right.
\tag {1.11}
$$
$a$ and $b$ are related to the original $\alpha$ and $\beta$ in Section 1.1 through the equations
$$\displaystyle
\left(
\begin{array}{c}
\displaystyle
a\\
\\
\displaystyle
b
\end{array}
\right)
=\frac{1}
{\sigma_0^2+\omega^2\epsilon_0^2}
\left(
\begin{array}{lr}
\displaystyle
\sigma_0 & \omega^2\epsilon_0\\
\\
\displaystyle
-\epsilon_0 & \sigma_0
\end{array}
\right)
\left(\begin{array}{c}
\displaystyle
\alpha\\
\\
\displaystyle
\beta
\end{array}
\right).
$$

\subsection{The enclosure method}

Let us recall notation and some definition.

We denote by $S^{n-1}$ the set of all unit vectors in $\Bbb R^n$.  The function $h_D$ defined by
the equation
$$\displaystyle
h_D(\vartheta)=\sup_{x\in D}x\cdot\vartheta,\,\,\vartheta\in S^{n-1}
$$
is called the support function of $D$.
For each $\vartheta\in S^{n-1}$ and a positive number $\delta$ set
$$\displaystyle
D_{\vartheta}(\delta)
=\{x\in D\,\vert\,h_D(\vartheta)-\delta<x\cdot\vartheta\le h_D(\vartheta)\}.
$$

{\bf\noindent Definition (Jump condition).}
Given $\vartheta\in S^{n-1}$ we say that
$\sigma$ has a positive jump on $\partial D$
from the direction $\vartheta$ if there exist constants $C_{\vartheta} > 0$ and $\delta_{\vartheta} > 0$
such that, for almost all $x\in D_{\vartheta}(\delta_{\vartheta})$ 
the lowest eigenvalue of $a(x)$ is greater than $C_{\vartheta}$;
$\sigma$ has a negative jump on 
$\partial D$ from the direction $\vartheta$ if there exist constants $C_{\vartheta} > 0$
and  $\delta_{\vartheta} > 0$ such that
for almost all $x\in D_{\vartheta}(\delta_{\vartheta})$  the lowest eigenvalue of 
$- a(x)$ is greater than $C_{\vartheta}$.
(1.12)

Given $\vartheta\in S^{n-1}$ take $\vartheta^{\perp}\in S^{n-1}$
perpendicular to $\vartheta$.  Given $\tau>0$ and $t\in\Bbb R$ define
$$\displaystyle
I_{\vartheta,\vartheta^{\perp}}(\tau,t)
=e^{-2\tau t}
\text{Re}\,
<(\Lambda_{\sigma,\epsilon}
-\Lambda_{1,0})
(e^{\tau x\cdot(\vartheta+i\vartheta^{\perp})}\vert_{\partial\Omega}),
\overline{e^{\tau x\cdot(\vartheta+i\vartheta^{\perp})}\vert_{\partial\Omega}}>.
\tag {1.13}
$$
In the theorems stated below we always assume that $\partial D$ is Lipschitz, $C^2$ in the case 
when $n=2, 3$, respectively.

\proclaim{\noindent Theorem 1.1.}
Assume that $\sigma$ has a positive jump on $\partial D$ from the direction $\vartheta$.
Then we have

 if $t>h_D(\vartheta)$, then 
$\displaystyle
\lim_{\tau\longrightarrow\infty}\vert I_{\vartheta,\vartheta^{\perp}}(\tau,t)\vert=0$;

 if $t<h_D(\vartheta)$, then 
$\displaystyle
\lim_{\tau\longrightarrow\infty}\vert I_{\vartheta,\vartheta^{\perp}}(\tau,t)\vert=\infty$;

 if $t=h_D(\vartheta)$, then 
$\displaystyle
\liminf_{\tau\longrightarrow\infty}\vert I_{\vartheta,\vartheta^{\perp}}(\tau,t)\vert>0$.

Moreover, the formula
$$\displaystyle
\lim_{\tau\longrightarrow\infty}
\frac{\log\vert I_{\vartheta,\vartheta^{\perp}}(\tau,t)\vert}
{2\tau}
=h_D(\vartheta)-t\,\,\,\,\,\,\forall t\in\Bbb R,
$$
is valid.

\endproclaim

Note that there is no restriction on $\omega$.
However, if $\sigma$ has a negative jump on $\partial D$ from
direction $\vartheta$, we do not know whether one can relax the condition (1.15) indicated below.

\proclaim{\noindent Theorem 1.2.}
Let $M>0$ and $m>0$ satisfy
$$\begin{array}{c}
\displaystyle
\sigma(x)\xi\cdot\xi\ge m\vert\xi\vert^2\,\,\text{a.e.}\,x\in D\,\,\,\,\forall\xi\in\Bbb R^n\\
\\
\displaystyle
\vert b(x)\xi\vert\le M\vert\xi\vert\,\,\text{a.e.}\,x\in D\,\,\,\,\forall\xi\in\Bbb R^n.
\end{array}
\tag {1.14}
$$
Assume that $\sigma$ has a negative jump on $\partial D$ from the direction $\vartheta$ and that, for the constant
$C_{\vartheta}$ in (1.12) the frequency $\omega$ satisfies
$$\displaystyle
0\le\omega<\frac{\sqrt{mC_{\vartheta}}}{M}.
\tag {1.15}
$$
Then we have the same conclusion as that of Theorem 1.1.

\endproclaim

In \cite{Ie2} we gave an extraction formula for the support function of polygonal inclusions in
the case when $n = 2$, $\omega=0$ and $a(x)$ is isotropic and constant from a single set of the Cauchy
data on $\partial\Omega$ of a solution of the governing equation. It would be interesting to considerwhether
the method still works or not. This remains open.

Algorithms for drawing a picture of the convex hull of $D$ based on formulae for $\omega=0$
are proposed in \cite{BH, IO, ISm1} and therein numerical testings are done. It would be interesting to
consider an algorithm for the purpose and do the numerical testing. The difference from the
previous situation is that one may use finitely many frequencies $\omega=\omega_1,\cdots,\omega_m$.

\subsection{A generalization}

In this section we employ the idea in [9] and consider the case when $n = 2$. Let $0<\alpha\le 1$.
The entire function
$$\displaystyle
E_{\alpha}(z)
=\sum_{n=0}^{\infty}\frac{z^n}{\Gamma(\alpha n+1)},\,\,z\in\mbox{\boldmath $C$}
$$
is called Mittag-Leffler's function (pages 206-208 on \cite{B} and \cite{ML}).
This includes $e^z$ as a special case because $E_1(z)=e^z$.  If $0<\alpha<1$, this function has the following remarkable property as $\vert z\vert\longrightarrow\infty$:

\noindent
if $\vert\text{arg}\,z\vert\le\pi\alpha/2$, then as $\vert z\vert\longrightarrow\infty$
$$\displaystyle
E_{\alpha}(z)\sim\frac{1}{\alpha}e^{z^{1/\alpha}};
$$

\noindent
if $\pi\alpha/2<\vert\text{arg}\,z\vert\le\pi$, then as $\vert z\vert\longrightarrow\infty$
$$\displaystyle
E_{\alpha}(z)\sim-\frac{z^{-1}}{\Gamma(1-\alpha)}.
$$

Let $y\in\Bbb R^2$ and $\vartheta\in S^1$.  Take $\vartheta^{\perp}\in S^1$ such that
$\vartheta\cdot\vartheta^{\perp}=0$.
For each $t\in\Bbb R$ consider the functions depending on $\tau>0$:
$$\displaystyle
e^{\alpha}_{\tau}(x;y,\vartheta,\vartheta^{\perp},t)
=E_{\alpha}(\tau\{(x-y)\cdot\vartheta-t+i(x-y)\cdot\vartheta^{\perp}\}).
$$
These functions are harmonic. 
Let $\mbox{$\cal C$}_{y+t\vartheta}(\vartheta,\pi\alpha/2)$ denote the cone about $\vartheta$ of 
opening angle $\pi\alpha/2$ with vertex at $y+t\vartheta$.

From the prperty of $E_{\alpha}(z)$ mentioned above one knows that

if $x\in\mbox{$\cal C$}_{y+t\vartheta}(\vartheta,\pi\alpha/2)
\setminus\partial\mbox{$\cal C$}_{y+t\vartheta}(\vartheta,\pi\alpha/2)$, then $\vert e^{\alpha}_{\tau}(x;y,\vartheta,\vartheta^{\perp},t)
\vert\longrightarrow\infty$ as $\tau\longrightarrow\infty$;

if $x\in\Bbb R^2\setminus\mbox{$\cal C$}_{y+t\vartheta}(\vartheta,\pi\alpha/2)$, then $\vert e^{\alpha}_{\tau}(x;y,\vartheta,\vartheta^{\perp},t)\vert\longrightarrow 0$ as $\tau\longrightarrow\infty$.

Define
$$\displaystyle
I^{\alpha}_{(y,\vartheta)}(\tau,t)
=\text{Re}\,
<(\Lambda_{\sigma,\epsilon}-\Lambda_{1,0})
(e^{\alpha}_{\tau}(x;y,\vartheta,\vartheta^{\perp},t)\vert_{\partial\Omega}),
\overline{e^{\alpha}_{\tau}(x;y,\vartheta,\vartheta^{\perp},t)\vert_{\partial\Omega}}>.
$$
Note that $I^{\alpha}_{(y,\vartheta)}(\tau,t)$ does not depend on the choice of $\vartheta^{\perp}$.

{\bf\noindent Definition (Generalized support function).}
Given $(y,\vartheta)\in (\Bbb R^2\setminus\overline\Omega)\times S^1$ with
$$\displaystyle
\mbox{$\cal C$}_{y}(\vartheta,\pi\alpha/2)\subset\Bbb R^2\setminus\Omega
\tag {1.16}
$$
define
$$\displaystyle
h^{\alpha}_D(y,\vartheta)
=\inf\left\{
t\in\,]-\infty,\,0[\,\vert\,
\forall s\in\,]t,\,0[\,
\mbox{$\cal C$}_{y+s\vartheta}(\vartheta,\pi\alpha/2)
\subset\Bbb R^2\setminus\overline D\right\}.
$$

The generalized support function gives the estimation of $D$ from above in the following sense:
$$\displaystyle
D\subset
\Bbb R^2\setminus
\overline{\mbox{$\cal C$}_{y+t\vartheta}(\vartheta,\pi\alpha/2)},\,\,t=h^{\alpha}_D(y,\vartheta).
$$
In the theorems stated below we always assume that $\partial D$ is Lipschitz.

\proclaim{\noindent Theorem 1.3.}
Let $0<\alpha<1$.  Assume that $a(x)$ is uniformly positive definite in $D$.
Let $(y,\vartheta)\in\,(\Bbb R^2\setminus\overline\Omega)\times S^1$ satisfy (1.16).

Then we have

if $t>h^{\alpha}_D(y,\vartheta)$, then $\lim_{\tau\longrightarrow\infty}\vert I^{\alpha}_{(y,\vartheta)}(\tau,t)\vert=0$;

if $t<h^{\alpha}_D(y,\vartheta)$, then $\lim_{\tau\longrightarrow\infty}\vert I^{\alpha}_{(y,\vartheta)}(\tau,t)\vert=\infty$;

if $t=h^{\alpha}_D(y,\vartheta)$, then $\liminf_{\tau\longrightarrow\infty}\vert I^{\alpha}_{(y,\vartheta)}(\tau,t)\vert>0$.

\endproclaim

This gives the characterization of $h^{\alpha}_D(y,\vartheta)$:
$$\displaystyle
]h^{\alpha}_D(y,\vartheta),\,0[
=\{t\in\,]-\infty,\,0[\,\vert\,
\lim_{\tau\longrightarrow\infty}I^{\alpha}_{(y,\vartheta)}(\tau,t)=0\}.
$$

\noindent
From Theorem 1.3 one knows that the set of all points on $\partial D$ that are visible from infinity 
can be reconstructed from $I^{\alpha}_{(y,\vartheta)}(\tau,t)$ for all $\alpha\in\,]0,\,1[$,
all $y$ on the circle with a large radius and all $\omega\in S^1$.  In the case when $\sigma_0=0$ in the original $\sigma$ of
Section 1.1, Theorem 3.1 says that one can extract the visible part of $\partial D$ from $\Lambda_{\sigma,\epsilon}$
provided, briefly speaking, the original $\epsilon$ of Section 1.1 is greater than $\epsilon_0$ in $D$.  There is no restriction on 
the bound of $\omega$.

\proclaim{\noindent Theorem 1.4.}
Let $0<\alpha<1$.  Let $M>0$ and $m>0$ satisfy (1.14).
Assume that there exists a positive number $C$ such that, for almost all $x\in D$, the lowest 
eigenvalue of $-a(x)$ is greater than $C$.  Let $(y,\vartheta)\in\,(\Bbb R^2\setminus\overline\Omega)\times S^1$ satisfy (1.16).
Let $\omega$ satisfy
$$\displaystyle
0\le\omega<\frac{\sqrt{mC}}{M}.
$$
Then we have the same conclusion as that of Theorem 1.3.

\endproclaim

We do not know whether one can relax the restriction on $\omega$.
This remains open.
The next problem is how to regularize the characterization of the generalized support function and
propose an algorithm based on the regularization for drawing a picture of giving an estimation of $D$ 
from above.  This remains open.

\subsection {Other related results}
In the case when $n = 2$, using a method in Brown-Uhlmann \cite{BU} and a perturbation argument,
Francini \cite{F} proved: if both $\sigma$ and $\epsilon$ are isotropic and have a regularity stronger than continuity,
then $\Lambda_{\sigma,\epsilon}$ uniquely determines $\sigma$ and $\epsilon$
themselves provided $\omega$ is small.

For drawing a picture of $D$ in the case when $\omega=0$ there is another interesting formula
established in \cite{BR}. It would be interesting to consider whether their method still works or not
for the case when $\omega\not=0$.

\section{A system of integral inequalities}

For the proof of theorems the system of integral inequalities (2.1) and (2.2) indicated below is
crucial.

\proclaim{\noindent Proposition 2.1.}  
Let $(\sigma_1,\epsilon_2)$ and $(\sigma_2,\epsilon_2)$ denote two pairs of conductivity and permittivity.
Assume that  both $\sigma_1$ and $\sigma_2$ are uniformly positive definite in $\Omega$.
Given $f\in H^{1/2}(\partial\Omega)$ let $u_j\in H^1(\Omega)$ denote the weak solution of
$$\begin{array}{c}
\displaystyle
\nabla\cdot(\sigma_j-i\omega\epsilon_j)\nabla u_j=0\,\,\text{in}\,\Omega,\\
\\
\displaystyle
u_j=f\,\,\text{on}\,\partial\Omega.
\end{array}
$$
Then we have
$$\begin{array}{c}
\displaystyle
\int_{\Omega}
(\sigma_1+i\omega\epsilon_1)\{(\sigma_1+\omega^2\epsilon_1\sigma_1^{-1}\epsilon_1)^{-1}
-\sigma_2^{-1}\}
(\sigma_1-i\omega\epsilon_1)\nabla u_1\cdot\overline{\nabla u_1}dx\\
\\
\displaystyle
\le
\text
{Re}\,<(\Lambda_{\sigma_2,\epsilon_2}-\Lambda_{\sigma_1,\epsilon_1})f,\overline f>;
\end{array}
\tag {2.1}
$$
$$\displaystyle
\text
{Re}\,<(\Lambda_{\sigma_2,\epsilon_2}-\Lambda_{\sigma_1,\epsilon_1})f,\overline f>
\le
\int_{\Omega}\{(\sigma_2+\omega^2\epsilon_2\sigma_2^{-1}\epsilon_2)-\sigma_1\}\nabla u_1\cdot\overline{\nabla u_1}dx.
\tag {2.2}
$$

\endproclaim

Note that, if $\omega=0$, then inequalities (2.1) and (2.2) coincide with those established in \cite{ISI}.

{\it\noindent Proof.}
Since $u_2=u_1=f$ on $\partial\Omega$ and $\overline{u_1}$ satisfies $\nabla\cdot(\sigma_1+i\omega\epsilon_1)\nabla\overline{u_1}
=0$ in $\Omega$,
we have
$$\begin{array}{c}
\displaystyle
\int_{\Omega}(\sigma_1-i\omega\epsilon_1)\nabla u_1\cdot\overline{\nabla(u_2-u_1)}dx=0;\\
\\
\displaystyle
\int_{\Omega}(\sigma_1-i\omega\epsilon_1)\nabla(u_2-u_1)\cdot\overline{\nabla u_1}dx=
-2i\omega\int_{\Omega}\epsilon_1\overline{\nabla u_1}\cdot\nabla(u_2-u_1)dx.
\end{array}
$$
Then, it is easy to see that
$$\begin{array}{c}
\displaystyle
<(\Lambda_{\sigma_2,\epsilon_2}-\Lambda_{\sigma_1,\epsilon_1})f,\overline f>
=\int_{\Omega}(\sigma_2-i\omega\epsilon_2)\nabla u_2\cdot\overline{\nabla u_2}dx
-\int_{\Omega}(\sigma_1-i\omega\epsilon_1)\nabla u_1\cdot\overline{\nabla u_1}dx
\\
\\
\displaystyle
=\int_{\Omega}
\{(\sigma_2-\sigma_1)-i\omega(\epsilon_2-\epsilon_1)\}
\nabla u_2\cdot\overline{\nabla u_2}dx
\\
\\
\displaystyle
+\int_{\Omega}(\sigma_1-i\omega\epsilon_1)\nabla u_2\cdot\overline{\nabla u_2}dx
-\int_{\Omega}(\sigma_1-i\omega\epsilon_1)\nabla u_1\cdot\overline{\nabla u_1}dx\\
\\
\displaystyle
=\int_{\Omega}\left\{(\sigma_1-i\omega\epsilon_1)\nabla(u_1-u_2)\cdot\overline{\nabla(u_1-u_2)}
+\{(\sigma_2-\sigma_1)-i\omega(\epsilon_2-\epsilon_1)\}\nabla u_2\cdot\overline{\nabla u_2}\right\}dx\\
\\
\displaystyle
+\int_{\Omega}(\sigma_1-i\omega\epsilon_1)\nabla u_1\cdot\overline{\nabla (u_2-u_1)}dx
+\int_{\Omega}(\sigma_1-i\omega\epsilon_1)\nabla (u_2-u_1)\cdot\overline{\nabla u_1}dx\\
\\
\displaystyle
=\int_{\Omega}\left\{(\sigma_1-i\omega\epsilon_1)\nabla(u_1-u_2)\cdot\overline{\nabla (u_1-u_2)}
+\{(\sigma_2-\sigma_1)-i\omega(\epsilon_2-\epsilon_1)\}
\nabla u_2\cdot\overline{\nabla u_2}\right\}dx\\
\\
\displaystyle
-2i\omega
\int_{\Omega}\epsilon_1\overline{\nabla u_1}\cdot\nabla(u_2-u_1)dx
\end{array}
\tag {2.3}
$$
and
$$\displaystyle
(\sigma_1-i\omega\epsilon_1)
\nabla(u_1-u_2)\cdot\overline{\nabla(u_1-u_2)}
+\{(\sigma_2-\sigma_1)-i\omega(\epsilon_2-\epsilon_1)\}
\nabla u_2\cdot\overline{\nabla u_2}
=A-i\omega B
$$
where
$$\begin{array}{c}
\displaystyle
A=\sigma_1\nabla(u_1-u_2)\cdot\overline{\nabla(u_1-u_2)}
+(\sigma_2-\sigma_1)\nabla u_2\cdot\overline{\nabla u_2};
\\
\\
\displaystyle
B=\epsilon_1\nabla(u_1-u_2)\cdot\overline{\nabla(u_1-u_2)}
+(\epsilon_2-\epsilon_1)\nabla u_2\cdot\overline{\nabla u_2}.
\end{array}
$$
Since $\sigma_j$, $\epsilon_j$ are real and symmetric, both $A$ and $B$ are real.
Therefore, we have
$$\displaystyle
\text{Re}\,<(\Lambda_{\sigma_2,\epsilon_2}
-\Lambda_{\sigma_1,\epsilon_1})f,\overline f>
=\int_{\Omega} Adx
-\int_{\Omega}
\text{Re}\,\{2i\omega\epsilon_1\overline{\nabla u_1}\cdot\nabla(u_2-u_1)\}dx.
$$
Write
$$\begin{array}{c}
\displaystyle
A-\text{Re}\,\{2i\omega\epsilon_1\overline{\nabla u_1}\cdot
\nabla(u_2-u_1)\}
=\sigma_2\nabla u_2\cdot\overline{\nabla u_2}
-(\sigma_1\nabla u_1\cdot\overline{\nabla u_2}+
\sigma\nabla u_2\cdot\overline{\nabla u_1})\\
\\
\displaystyle
+\sigma_1\nabla u_1\cdot\overline{\nabla u_1}
-(i\omega\epsilon_1\overline{\nabla u_1}\cdot\nabla u_2-i\omega\epsilon_1\nabla u_1\cdot\overline{\nabla u_2})\\
\\
\displaystyle
=\sigma_2\nabla u_2\cdot\overline{\nabla u_2}
-(\sigma_1-i\omega\epsilon_1)\nabla u_1\cdot\overline{\nabla u_2}
-(\sigma_1+i\omega\epsilon_1)\overline{\nabla u_1}\cdot\nabla u_2
+\sigma_1\nabla u_1\cdot\overline{\nabla u_1}.
\end{array}
$$
Then we have two inequalities:
$$\begin{array}{c}
\displaystyle
A-\text{Re}\,\{2i\omega\epsilon_1\overline{\nabla u_1}\cdot\nabla(u_2-u_1)\}
=\vert\sigma_2^{1/2}\nabla u_2-\sigma_2^{-1/2}(\sigma_1-i\omega\epsilon_1)\nabla u_1\vert^2\\
\\
\displaystyle
+\sigma_1\nabla u_1\cdot\overline{\nabla u_1}-\vert\sigma_2^{-1/2}(\sigma_1-i\omega\epsilon_1)\nabla u_1\vert^2
\\
\\
\displaystyle
\ge
\{\sigma_1-(\sigma_1+i\omega\epsilon_1)\sigma_2^{-1}(\sigma_1-i\omega\epsilon_1)\}\nabla u_1\cdot\overline{\nabla u_1};
\end{array}
\tag {2.4}
$$
$$\begin{array}{c}
\displaystyle
A-\text{Re}\,\{2i\omega\epsilon_1\overline{\nabla u_1}\cdot\nabla(u_2-u_1)\}
=\vert\sigma_1^{1/2}\nabla u_1-\sigma_1^{-1/2}(\sigma_1+i\omega\epsilon_1)\nabla u_2\vert^2\\
\\
\displaystyle
+\sigma_2\nabla u_2\cdot\overline{\nabla u_2}-\vert\sigma_1^{-1/2}(\sigma_1+i\omega\epsilon_1)\nabla u_2\vert^2
\\
\\
\displaystyle
\ge
\{\sigma_2-(\sigma_1-i\omega\epsilon_1)\sigma_1^{-1}(\sigma_1+i\omega\epsilon_1)\}\nabla u_2\cdot\overline{\nabla u_2}.
\end{array}
\tag {2.5}
$$
Then from (2.4) and the identity
$$\begin{array}{c}
\displaystyle
\sigma_1-(\sigma_1+i\omega\epsilon_1)\sigma_2^{-1}(\sigma_1-i\omega\epsilon_1)\\
\\
\displaystyle
=(\sigma_1+i\omega\epsilon_1)\{(\sigma_1+i\omega\epsilon_1)^{-1}\sigma_1(\sigma_1-i\omega\epsilon_1)^{-1}
-\sigma_2^{-1}\}(\sigma_1-i\omega\epsilon_1)\\
\\
\displaystyle
=(\sigma_1+i\omega\epsilon_1)\{((\sigma_1-i\omega\epsilon_1)\sigma_1^{-1}(\sigma_1+i\omega\epsilon_1))^{-1}
-\sigma_2^{-1}\}(\sigma_1-i\omega\epsilon_1)\\
\\
\displaystyle
=(\sigma_1+i\omega\epsilon_1)\{(\sigma_1+\omega^2\epsilon_1\sigma^{-1}\epsilon_1)^{-1}-\sigma_2^{-1}\}
(\sigma_1-i\omega\epsilon_1)
\end{array}
$$
we obtain (2.1).  From (2.5) we obtain
$$\begin{array}{c}
\displaystyle
\text{Re}\,<(\Lambda_{\sigma_2,\epsilon_2}-\Lambda_{\sigma_1,\epsilon_1})f,\overline f>
\ge
\int_{\Omega}
\{(\sigma_2-(\sigma_1-i\omega\epsilon_1)\sigma_1^{-1}
(\sigma_1+i\omega\epsilon_1)\}
\nabla u_2\cdot\overline{\nabla u_2}dx.
\end{array}
\tag {2.6}
$$
By interchanging subscripts $1$ and $2$, we obtain (2.2).

\noindent
$\Box$

{\bf\noindent Remark 2.1.}
It is easy to see that, for any invertible matrices $A$ and $B$ we have the identity
$$\displaystyle
A^{-1}-B^{-1}
=B^{-1}(B-A)B^{-1}
+B^{-1}(B-A)A^{-1}(B-A)B^{-1}.
$$
This yields
$$\begin{array}{c}
\displaystyle
(\sigma_1+\omega^2\epsilon_1\sigma_1^{-1}\epsilon_1)^{-1}
-\sigma_2^{-1}
=\sigma_2^{-1}\{\sigma_2-(\sigma_1+\omega^2\epsilon_1\sigma_1^{-1}\epsilon_1)\}\sigma_2^{-1}\\
\\
\displaystyle
+\sigma_2^{-1}
\{\sigma_2-(\sigma_1+\omega^2\epsilon_1\sigma_1^{-1}\epsilon_1)\}
(\sigma_1+\omega^2\epsilon_1\sigma_1^{-1}\epsilon_1)^{-1}
\{\sigma_2-(\sigma_1+\omega^2\epsilon_1\sigma_1^{-1}\epsilon_1)\}\sigma_2^{-1}.
\end{array}
\tag {2.7}
$$

\section{Proof of theorems.}

Let $v\in H^1(\Omega)$ satisfy
$$\begin{array}{c}
\displaystyle
\nabla\cdot\nabla v=0\,\,\text{in}\,\Omega,\\
\\
\displaystyle
v=f\,\,\text{on}\,\partial\Omega.
\end{array}
$$
Set $(\sigma_1,\epsilon_1)=(1,0)$ and $(\sigma_2,\epsilon_2)=(\sigma,\epsilon)$.
Then from (2.1) and (2.2) we obtain
$$\displaystyle
\int_D(1-\sigma^{-1})\nabla v\cdot\overline{\nabla v}dx
\le
\text{Re}\,<(\Lambda_{\sigma,\epsilon}-\Lambda_{1,0})f,\overline f>;
\tag {3.1}
$$
$$\displaystyle
\text{Re}\,<(\Lambda_{\sigma,\epsilon}-\Lambda_{1,0})f,\overline f>
\le\int_D(\sigma+\omega^2\epsilon\sigma^{-1}\epsilon-1)\nabla v\cdot\overline{\nabla v}dx.
\tag {3.2}
$$

{\it\noindent Proof of Theorems 1.1 and 1.2.}
First we consider the case:
$$\displaystyle
\text{for almost all $x\in D_{\vartheta}(\delta_{\vartheta})$ the lowest eigenvalue of $-a(x)$ is greater than $C_{\vartheta}$.}
\tag {3.3}
$$
From (3.1) and (3.2) for $v=e^{\tau x\cdot(\vartheta+i\vartheta^{\perp})}$, one can easily see that
$$\displaystyle
I_{\vartheta,\vartheta^{\perp}}(\tau,t)\vert_{t=h_D(\vartheta)}
\le e^{-2\tau h_D(\vartheta)}
\int_D(\sigma+\omega^{2}\epsilon\sigma^{-1}\epsilon-1)\nabla v\cdot\overline{\nabla v}dx;
\tag {3.4}
$$
$$\displaystyle
I_{\vartheta,\vartheta^{\perp}}(\tau,t)\vert_{t=h_D(\vartheta)}
=O(\tau^2)
\tag {3.5}
$$
as $\tau\longrightarrow\infty$.

From the assumption on the regularity for $\partial D$, one can find constants $M_{\vartheta}>0$,
$\epsilon_{\vartheta}>0$
such that
$$\displaystyle
\mu_{n-1}(\{x\in D\,\vert\,x\cdot\vartheta=h_D(\vartheta)-s\})\ge M_{\vartheta}s,\,\,\forall s\in\,]0,\,\epsilon_{\vartheta}[
\tag {3.6}
$$
where $\mu_{n-1}$ denotes the $(n-1)$-dimensional Lebesgue measure.
Set $\delta'\equiv\delta'_{\vartheta}
=\min \{\delta_{\vartheta},\epsilon_{\vartheta}\}$.
Since $D'\equiv D_{\vartheta}(\delta')\subset D_{\vartheta}(\delta_{\vartheta})$,
from (1.14) we have
$$\begin{array}{c}
\displaystyle
\int_{D'}(\sigma+\omega^2\epsilon\sigma^{-1}\epsilon-1)\nabla v\cdot\overline{\nabla v}dx
=\int_{D'}(a(x)\nabla v\cdot\overline{\nabla v}+
\omega^2\sigma(x)^{-1}b(x)\nabla v\cdot\overline{b(x)\nabla v}dx\\
\\
\displaystyle
\le
-(C_{\vartheta}-\omega^2m^{-1}M^2)
\int_{D'}\vert\nabla v\vert^2dx
=-2(C_{\vartheta}-\omega^2m^{-1}M^2)\tau^2\int_{D'}e^{2\tau x\cdot\vartheta} dx.
\end{array}
\tag {3.7}
$$
Using (3.6), we obtain
$$\begin{array}{c}
\displaystyle
e^{-2\tau h_D(\vartheta)}\tau^2
\int_{D'}e^{2\tau x\cdot\vartheta}dx
=\tau^2\int_0^{\delta'}e^{-2\tau s}
\mu_{n-1}(\{x\in D\,\vert x\cdot\vartheta=h_D(\vartheta)-s\})ds\\
\\
\displaystyle
\ge M_{\vartheta}\tau^2
\int_0^{\delta'}se^{-2\tau s}ds
=M_{\vartheta}\int_0^{\tau\delta'}te^{-2t}dt.
\end{array}
\tag {3.8}
$$
A combination of (3.7) and (3.8) yields, given $\tau_0>0$
$$\displaystyle
e^{-2\tau h_D(\vartheta)}
\int_{D'}
(\sigma+\omega^2\epsilon\sigma^{-1}\epsilon-1)\nabla v\cdot\overline{\nabla v}dx
\le -2(C_{\vartheta}-\omega^2m^{-1}M^2)C_{\tau_0\delta'}M_{\vartheta}
\tag {3.9}
$$
for all $\tau\ge\tau_0$.  Since $x\cdot\vartheta\le h_D(\vartheta)-\delta'$ for all $x\in D\setminus D'$, it is easy to see that
$$\displaystyle
e^{-2\tau h_D(\vartheta)}
\int_{D\setminus D'}(\sigma+\omega^2\epsilon\sigma^{-1}\epsilon-1)\nabla v\cdot\overline{\nabla v}dx=O(\tau^2e^{-2\tau\delta'}).
\tag {3.10}
$$
From (3.9) and (3.10) one concludes that if $\tau_0$ is sufficiently large, then for all $\tau\ge\tau_0$ one has
$$\displaystyle
e^{-2\tau h_D(\vartheta)}\int_D(\sigma+\omega^2\epsilon\sigma^{-1}\epsilon-1)\nabla v\cdot\overline{\nabla v}dx
\le -C
\tag {3.11}
$$
where $C$ is a positive constant provided $\omega$ satisfies (1.15).

Now everything comes from (3.4), (3.5), (3.11) and the trivial identity
$$\displaystyle
I_{\vartheta,\vartheta^{\perp}}(\tau,t)
=e^{2\tau(h_D(\vartheta)-t)}I_{\vartheta,\vartheta^{\perp}}(\tau,t)\vert_{t=h_D(\vartheta)}.
$$
Using (3.1) and (2.7) for $(\sigma_1,\epsilon_1)=(1,0)$ and $(\sigma_2,\epsilon_2)=(\sigma,\epsilon)$,
one can also complete the proof in the case when
$$\displaystyle
\text{for almost all $x\in D_{\vartheta}(\delta_{\vartheta})$ the lowest eigenvalue of $a(x)$ is greater than $C_{\vartheta}$.}
$$
Note that in this case we do not make use of (1.15).

\noindent
$\Box$

{\it\noindent Proof of Theorems 1.3 and 1.4.}
We briefly describe the outline of the proof.
Using inequalities (3.1), (3.2) for $v=e^{\alpha}_{\tau}(x;y,\vartheta,\vartheta^{\perp},t)$, we obtain
$$\displaystyle
C_1J(\tau,t)\le\vert I_{(y,\vartheta)}^{\alpha}(\tau,t)\vert
\le C_2J(\tau,t)
$$
where $C_1$, $C_2$ are positive constants independent of $y,\vartheta,\vartheta^{\perp},t,\tau$;
$$\displaystyle
J(\tau,t)
=\int_D\vert\nabla e^{\alpha}_{\tau}(x;y,\vartheta,\vartheta^{\perp},t)\vert^2dx
=2\tau^2\int_D\vert E'_{\alpha}(\tau\{(x-y)\cdot\vartheta-t+i(x-y)\cdot\vartheta^{\perp}\})\vert^2dx.
$$
Therefore, everything comes from the facts that

if $t>h^{\alpha}_D(y,\vartheta)$, then $\lim_{\tau\longrightarrow\infty}J(\tau,t)=0$;

if $t<h^{\alpha}_D(y,\vartheta)$, then $\lim_{\tau\longrightarrow\infty}J(\tau,t)=\infty$;

if $t=h^{\alpha}_D(y,\vartheta)$, then $\liminf_{\tau\longrightarrow\infty}J(\tau,t)>0$.

These facts are proved in \cite{IM} and we omit the proof.

\noindent
$\Box$

$$\quad$$

\centerline{{\bf Acknowledgments}}

The author started this research during the visit at
Mathematical Sciences Research Institute,
Berkeley in August, 2001.
This research was partially supported by Grant-in-Aid for
Scientific Research (C)(No. 13640152) of Japan  Society for
the Promotion of Science.

The author would like to thank the referees and the board member for several suggestions.

$$\quad$$

\end{document}